\documentclass{article}
\usepackage{amsfonts,amsbsy, amsmath, amssymb, amsthm}
\usepackage{verbatim}

\title{Cohomological invariants of odd degree Jordan algebras}

\author{Mark L. MacDonald}

\begin{document}

\newtheorem{thm}{Theorem}[section]
\newtheorem{lem}[thm]{Lemma}
\newtheorem{prop}[thm]{Proposition}
\newtheorem{cor}[thm]{Corollary}
\newtheorem{claim}[thm]{Claim}
\theoremstyle{remark}
\newtheorem{remark}[thm]{\textbf{Remark}}
\newtheorem{ex}[thm]{\textbf{Example}}

\maketitle

\begin{abstract}
In this paper we determine all possible cohomological invariants
of $\mathbf{Aut}(J)$-torsors in Galois cohomology with mod 2
coefficients (characteristic of the base field not 2), for $J$ a
split central simple Jordan algebra of odd degree $n\geq 3$. This
has already been done for $J$ of orthogonal and exceptional type,
and we extend these results to unitary and symplectic type. We
will use our results to compute the essential dimensions of some
groups, for example we show that ed$(\mathbf{PSp}_{2n})=n+1$ for
$n$ odd.
\end{abstract}

\section{Introduction}

The Stiefel-Whitney classes of quadratic forms over $k$ define
invariants in Galois cohomology $H^*(k, \mathbb{Z}/2\mathbb{Z})$
up to isometry \cite{De62}, \cite{Mi70}. It is shown in \cite[ch.\
VI]{GMS03} that the even Stiefel-Whitney classes form a basis of
all cohomological invariants of $\mathbf{SO}_n$-torsors for $n\geq
3$ odd. This is done by identifying the $\mathbf{SO}_n$-torsors
with isomorphism classes of determinant 1 quadratic forms of
dimension $n$. These torsors may be further identified with
isomorphism classes of algebras with orthogonal involution of
degree $n$, by sending $q$ to its adjoint involution $ad_q$ (see
\cite[Thm.\ 4.2, p.42]{KMRT}). These classes may further be
identified with isomorphism classes of central simple Jordan
algebras of degree $n$ whose associated composition algebra $C$ is
one dimensional (see \cite[p.210]{Ja68}). We wish to extend these
results to dim$(C)=2$ or 4, which is to say, odd degree algebras
with unitary and symplectic involutions. In fact, in the octonion
case dim$(C)=8$, we only have a Jordan algebra when the degree is
3, and then it is called an Albert algebra. The mod 2
cohomological invariants of Albert algebras have been determined
in \cite[ch.\ VI]{GMS03}. Nevertheless, we include this case here
for completeness.

For any $n\geq 3$ we will define $J^r_n:=(M_n(C),-)_+$ over $k$,
as the split Jordan algebra of hermitian elements, where $C$ is
the split composition algebra of dim$(C)=2^r$. If $r=3$ then we
insist that $n=3$. Here $-$ denotes the conjugate transpose
involution. The following table summarizes, for $n=2m+1$, the
(split) automorphism groups $\mathbf{Aut}(J^r_n)$, together with
their mod 2 cohomological invariants (see Thm.\ \ref{AutInvs}). We
list the degrees of invariants which form an $H^*(k_0)$-basis of
all invariants.

$$ $$

\begin{tabular}{c|c|c}
  $r$ & $\mathbf{Aut}(J^r_{2m+1})$ & Degrees of $H^*$-basis of Inv$(\mathbf{Aut}(J^r_n))$ \\
    \hline
 0 & $\mathbf{SO}_{2m+1}$ & 0, 2, 4, $\cdots$ , $2m$ \\
 1 & $\mathbb{Z}/2\mathbb{Z} \ltimes \mathbf{PGL}_{2m+1}$ & 0, 1, 3, $\cdots$ , $2m+1$ \\
 2 & $\mathbf{PSp}_{2(2m+1)}$ & 0, 2, 4, $\cdots$ , $2m+2$ \\
 3 & $F_4$ & 0, 3, 5 \\
\end{tabular}

$$ $$

Here $m \geq 1$, and $F_4$ denotes the split simple group of type
$F_4$. To specify the semi-direct product in $r=1$, we just need
to describe how the non-trivial element of
$\mathbb{Z}/2\mathbb{Z}$ acts on $\mathbf{PGL}_{2m+1}$. This
action is defined by sending any $[a]\in \mathbf{PGL}_{2m+1}$ to
its inverse transpose $[(a^t)^{-1}]$ (see \cite[29.20]{KMRT}).

We will show (Thm.\ \ref{AutInvs}) that for $r=1,2$ and 3 and
$n\geq 3$ odd, the normalized invariants of $\mathbf{Aut}(J^r_n)$
in degree greater than zero are in one-to-one correspondence with
the even Stiefel-Whitney classes of $n$-dimensional quadratic
forms (which are the invariants of $\mathbf{Aut}(J^0_n)$). Under
this bijection the degree zero Stiefel-Whitney class corresponds
to the degree $r$ invariant which classifies the associated
composition algebra.

The $\mathbf{Aut}(J^1_3)$-invariants are discussed in detail in
\cite[\S19.B, \S30.C]{KMRT} or \cite{HKRT}, including a mod 3
invariant of degree 2. These invariants determine an
$\mathbf{Aut}(J^1_3)$-torsor up to isomorphism (considered there
as a degree 3 algebra with unitary involution).

Recently in \cite{GPT07} Garibaldi, Parimala and Tignol have
classified mod 2 invariants of degree $\leq 3$ for
$\mathbf{Aut}(J^2_n)$-torsors for $n$ even.

In the final section we determine the essential dimension at 2 for
our groups $\mathbf{Aut}(J^r_n)$ for $n\geq 3$ odd. In each case
it is equal to the lower bound given by \cite[Theorem 1]{CS06}. In
particular, we find that ed$(\mathbf{PSp}_{2n})=n+1$ for $n\geq 3$
odd, where previously the best upper bound was given by
$2n^2-3n-1$ in \cite{Le04}.

\section{Preliminaries}

Throughout we will let $k$ be a field extension of a fixed base
field $k_0$ of characteristic not 2, and $k_s$ will be a separable
closure of $k$.

For any $k$-algebra $A$ (by which we will mean finite dimensional,
not necessarily associative, with identity), we will use the usual
notions from Galois cohomology \cite[\S29]{KMRT} to identify
$\mathbf{Aut}_{alg}(A)$-torsors over $k$ with $k$-isomorphism
classes of $k$-algebras $B$ such that $A_{k_s}\cong B_{k_s}$.
Similarly for algebras with involution.

In the Introduction we defined the split Jordan algebras
$J^r_n:=(M_n(C),-)_+$ for $r=0,1,2$ when $n\geq 3$, and also for
$r=3$ when $n=3$. Here, and throughout this paper we will write
dim$(C)=2^r$. These Jordan algebras are pairwise non-isomorphic,
and over $k_s$ they represent nearly all simple Jordan algebras by
a theorem of Albert (see \cite[Ch.\ V.6, p.204]{Ja68}). We will
say that a simple Jordan algebra is of degree $n$ (for some $n\geq
3$), if it becomes isomorphic to $J^r_n$ over $k_s$. These are the
only kind of Jordan algebras that we will consider.

Furthermore, for $r=0,1,2$ we have that $(M_n(C),-)$ is a central
simple algebra with orthogonal (resp.\ unitary, symplectic)
involution, where $C$ is again the split composition algebra of
dimension $2^r$. They are pairwise non-isomorphic, and over $k_s$
they form all $k_s$-isomorphism classes of central simple algebras
with involution. We say that a central simple algebra with
involution has degree $n$ if it becomes isomorphic to $(M_n(C),-)$
over $k_s$.

One must be careful with the potentially confusing terminology
here. A central simple algebra with involution over $k$ is defined
to be central and simple as an algebra-with-involution, and might
not be central or simple as an algebra over $k$ (see
\cite[p.208]{Ja68} or \cite[\S2]{KMRT} for precise definitions).

For $r=0,1,2$ we have a 1-to-1 correspondence between isomorphism
classes of these two types of objects given by $(A,\sigma)
\leftrightarrow (A,\sigma)_+$. So we can view
$\mathbf{Aut}(J^r_n)$-torsors as either isomorphism classes of
algebras with involution or as isomorphism classes of Jordan
algebras (see \cite[Ch.\ V.7, p.210]{Ja68}). An advantage of the
Jordan algebra point of view is that it includes the exceptional
$r=3$ case.

\subsection{Invariants of quadratic forms}

We will follow the notation of \cite{GMS03} and write $H^*(k)$ or
even $H^*$ for the Galois cohomology ring
$H^*(\textup{Gal}(k_s/k),\mathbb{Z}/2\mathbb{Z})$. For $a\in
k^*/(k^*)^2$, we will denote the corresponding element $(a)\in
H^1(k)$ so that we have $(a\cdot b)=(a)+(b)$. We let
$\mathsf{Quad}_{n,1}(k)$ be the pointed set of $k$-isometry
classes of $n$-dimensional quadratic forms of determinant 1. And
we let $\mathsf{Pf_r}(k)$ be the pointed set of $k$-isometry
classes of $r$-Pfister forms. We will write
$\textup{Inv}(G)=\textup{Inv}_{k_0}(H^1(-,G))$ for the group of
cohomological invariants in $H^*$ of $G$-torsors.

To define the Stiefel-Whitney classes of a quadratic form $q$ over
$k$, take a diagonalization $q \cong \langle a_1, \cdots , a_n
\rangle$. Then define $w_i(q)$ to be the $i$th degree part of the
product $$w(q)=\prod_{j=1}^n (1+(a_j)) \in H^*(k).$$ This product
is called the \textit{total Stiefel-Whitney class}. It is
independent of the diagonalization, which can be shown by a chain
equivalence argument (see \cite{Mi70}).

For an $r$-Pfister form, $$q=\langle \langle a_1, \cdots , a_r
\rangle \rangle := \langle 1, -a_1 \rangle \otimes \cdots \otimes
\langle 1, -a_r \rangle,$$ define an invariant by
$e_r(q)=(a_1)\cdots (a_r) \in H^r(k)$. It is shown in
\cite[VI]{GMS03} that the $H^*$-module of invariants of
$r$-Pfister forms, $\textup{Inv}(\mathsf{Pf}_r)$, has an
$H^*(k_0)$-basis consisting of $\{1, e_r\}$.

\section{An upper bound for the invariants}

In this section we show that Inv$(\mathbf{Aut}(J^r_n))$
injectively embeds into the tensor product of two groups of
invariants that we understand well (see Cor.\ \ref{InvInject}).

For any associative composition algebra $C$ over $k$, and
hermitian form $h$ on a free $n$-dimensional $C$-module $V$, the
\textit{trace form}, $q_h$ is the quadratic form over $k$ on $V$
such that $q_h(x)=h(x,x)$. If $\phi$ is the norm of $C$, and $V
\cong C \otimes V_0$ as $k$-vector spaces, then $q_h \cong \phi
\otimes q_0$ for some $n$-dimensional quadratic form $q_0$ on
$V_0$.

\begin{lem} \label{h-q_h}
Let $C$ be an associative composition algebra over $k$ with
conjugation involution, and let $h$, $h'$ be hermitian forms over
$C$. Then $h \cong h'$ iff $q_h \cong q_{h'}$.
\end{lem}

\begin{proof}
This is shown in \cite[10.1.1,10.1.7]{Sch85}.
\end{proof}

\begin{prop} \label{ad_h-h}
Let $C$ be an associative composition algebra over $k$ of
dimension $1$ (resp.\ $2$, $4$). Then isomorphism classes of
involutions on $M_n(C)$ of orthogonal (resp.\ unitary, symplectic)
type correspond to similarity classes of $n$-dimensional hermitian
forms over $C$, under $ad_h \leftrightarrow h$.
\end{prop}

\begin{proof}
This is proved in \cite[Prop.\ 4.2, p.43]{KMRT}.
\end{proof}

We will call a Jordan algebra over $k$ \textit{reduced} if it is
isomorphic to one of the form $(M_n \otimes C,ad_q \otimes -)_+$
for some composition $k$-algebra $C$, and $n$-dimensional
quadratic form $q$ over $k$.

Notice this implies Jacobson's \cite{Ja68} definition of reduced
in terms of orthogonal idempotents, and his definition implies
this one, by his Coordinatization theorem for $n\geq 3$.

\begin{lem} \label{ReducedAlg}
Let $J$ be a central simple Jordan algebra of odd degree $n$ and
of type $r=0$ or $2$. Then $J$ is reduced.
\end{lem}

\begin{proof}
So $J\cong (A,\sigma)_+$, where $(A,\sigma)$ is a central simple
algebra with an involution of the first kind. So its index is a
power of 2 (see \cite[p.18, 2.8]{KMRT}).

\textit{Case $r=0$}: The degree of $A$ as a central simple algebra
is $n$, so the index divides $n$, and hence $A$ is split. So $J$
is reduced by Prop.\ \ref{ad_h-h} (or \cite[p.1]{KMRT}).

\textit{Case $r=2$}: The degree of $A$ as a central simple algebra
is $2n$, so the index divides 2, and hence $J$ is reduced by
Prop.\ \ref{ad_h-h} and the remarks at the beginning of this
section.
\end{proof}

\begin{lem} \label{ReducingFields}
Let $J$ be a central simple Jordan algebra of odd degree. Then $J$
becomes reduced after extending scalars by a field extension of
odd degree.
\end{lem}

\begin{proof}
The only non-reduced algebras are when $r=3$ or $r=1$ by Lemma
\ref{ReducedAlg}. For the $r=3$ case, as stated in the
introduction, we must have $n=3$. By \cite[6.1]{SV00} any
non-reduced Albert algebra becomes reduced after a cubic
extension.

For $r=1$ there are two cases, as follows.

\textit{Case $J \cong (B \times B^{op}, \sigma)_+$}: Here $\sigma$
is the exchange involution, and $B$ is a central simple algebra
over $k$ of degree $n$ odd. Then any maximal subfield $L$ is a
splitting field for $B$ of degree dividing $n$. Then $J_L$ is
reduced by Prop.\ \ref{ad_h-h}.

\textit{Case $J \cong (A, \sigma)_+$}: Here $A$ is a central
simple algebra over $K$ of degree $n$ odd, where $K$ is a
quadratic extension of $k$, and $\sigma$ is an involution of the
second kind (i.e. unitary) over $k$. Since the Brauer group of a
finite field is trivial, we may assume that $k$ is infinite.

Let $L/k$ be an odd degree extension such that ind$(A_L)=d$ is
minimal, where $A_L$ is a simple associative algebra with centre
$K_L=K\otimes_k L$. $d$ is odd since $d|n$. Let $D$ be the
division $K_L$-algebra that is Brauer equivalent to $A_L$, and
hence of degree $d$ as a central simple algebra. Then $D$ has an
involution of the second kind $\tau$ that fixes $L$ by \cite[3.1,
p.31]{KMRT}.

We want to show that $d=1$, so assume that $d>1$. Then we can take
a non-scalar element $a$ in the degree $d$ Jordan algebra
$(D,\tau)_+$. Since $k$ is infinite, we may choose $a$ to be of
maximal degree in the sense of \cite[p.224]{Ja68}. In other words,
$a$ is such that deg$(m_a)=d$, where $m_a(\lambda)\in L[\lambda]$
is the minimal polynomial of $a$ in $(D,\tau)_+$, for some
indeterminate $\lambda$. Here we are using the fact from
\cite[p.233]{Ja68} that the degree of the generic minimal
polynomial of a generic element is equal to the degree of the
Jordan algebra as defined in the Preliminaries. Also, the
coefficients of $m_a$ are in $L$ by \cite[32.1.2, p.452]{KMRT}.

If we let $\alpha$ be a root of $m_a$ in $k_s$, then the minimal
polynomial of $\alpha$ is $m_a$. This is because $D$ is a division
algebra and hence $m_a$ is irreducible. So $E=L(\alpha)$ is a
field extension of degree $d$ over $L$, and in particular, of odd
degree over $k$. Then by considering the generic norm of
$a':=\alpha 1 - a_E \in D_E$, we see that $n(a')=m_a(\alpha)=0$
\cite[p.224]{Ja68}. So $a'\neq 0$ and is non-invertible, and hence
$D_E$ is not a division algebra. So
ind$(A_E)=\textup{ind}(D_E)<d$, which contradicts the minimality
of $d$.

Therefore $d=1$, and hence $A_L\cong M_n(K_L)$. So by Prop.
\ref{ad_h-h} we see that $J_L$ is reduced.
\end{proof}

\begin{prop} \label{Surject}
For $n\geq 3$ odd, let $J$ be an $\mathbf{Aut}(J^r_n)$-torsor over
$k$. Then there is an odd-degree extension $L/k$ such that $J_L$
is in the image of \begin{align*} \mathcal{H}: \mathsf{Pf}_r(L)
\times \mathsf{Quad}_{n,1}(L) &\to H^1(L,\mathbf{Aut}(J^r_n))\\
(\phi, q)\textup{    } &\mapsto (M_n \otimes C_{\phi}, ad_q
\otimes -)_+.
\end{align*} Moreover, if $r=2$ the map is a surjection, and for
$r=0$ it is a bijection.
\end{prop}

\begin{proof}
From Lemma \ref{ReducedAlg} and \ref{ReducingFields} we get $L/k$
such that $J_L$ is reduced, and since $n$ is odd, we can scale $q$
so that $det(q)=1$. Lemma \ref{ReducedAlg} gives the $r=2$
surjection, and the $r=0$ bijection is well-known.
\end{proof}

\begin{cor} \label{InvInject}
We have an injective map of invariants
$$\textup{Inv}(\mathbf{Aut}(J^r_n)) \hookrightarrow
\textup{Inv}(\mathsf{Pf}_r) \otimes
\textup{Inv}(\mathsf{Quad}_{n,1}).$$
\end{cor}

\begin{proof}
By \cite[Lemma 5.3]{Ga06} we can use the surjectivities from
Prop.\ \ref{Surject} to induce an injective map on invariants.
Then from \cite[Ex. 16.5]{GMS03}, we can express the invariants of
the direct product $\mathsf{Pf}_r \times \mathsf{Quad}_{n,1}$ as
the tensor product of the invariants of each factor.
\end{proof}

\section{Construction of the invariants}

Now it is a matter of deciding which of these invariants occur. In
other words, we wish to determine the image of the injective map
in Cor.\ \ref{InvInject}. It turns out to be the constant
invariants together with all multiples of $e_r$, the degree $r$
invariant of $\mathsf{Pf}_r$ (Thm.\ \ref{AutInvs}).

\begin{thm} \label{e_r.w}
For $n$ odd, $r=0,1,2$ or $3$, the invariants $e_r \otimes w_{2i}
\in \textup{Inv}(\mathsf{Pf}_r) \otimes
\textup{Inv}(\mathsf{Quad}_{n,1})$ extend uniquely to
$\mathbf{Aut}(J^r_n)$-invariants of degree $r+2i$, which we will
call $v_i$.
\end{thm}

If the invariants extend, then by Cor.\ \ref{InvInject} they are
unique. First we will show how to construct the invariants on
reduced Jordan algebras, and then use \cite[Prop.\ 7.2]{Ga06} to
extend them to all Jordan algebras.

For a reduced Jordan algebra $J=(M_n(C),ad_q \otimes -)_+$, we
call $C$ the \textit{composition algebra associated to $J$}. It is
determined up to isomorphism by the isomorphism class of $J$ (see
\cite{Ja68} or \cite[16]{Mc04}). We will usually denote its norm
form $\phi$, which is a Pfister form.

\begin{lem} \label{UptoSimilarity}
The quadratic form $\phi \otimes q$ is determined up to similarity
by the isomorphism class of the reduced Jordan algebra
$J=(M_n(C),ad_q \otimes -)_+$.
\end{lem}

\begin{proof}
First consider the case when $C$ is associative, which is to say,
$r \neq 3$. By \cite[p.210]{Ja68} the Jordan algebra $J$
determines the isomorphism class of the algebra with involution
$(M_n(C),\sigma)$. Then Prop.\ \ref{ad_h-h} lets us associate up
to a scalar, the $n$-dimensional hermitian form $h$ on $C$.
Finally, Lemma \ref{h-q_h} allows us to determine its trace form
$\phi \otimes q$ up to similarity.

In the case $r=3$, we use the following argument. For any reduced
Jordan algebra of degree $n$, the quadratic (reduced) trace form,
$T_J(x)=Trd(x^2)=trace(x^2)$ is determined up to isometry by the
isomorphism class of $J$, and is of the form $$T_J \cong n\langle
1 \rangle \perp \langle 2 \rangle \phi \otimes \wedge^2q.$$ But
since $q$ is similar to $\wedge^2(q)$ for $3$-dimensional forms,
we see that $\phi \otimes q$ is determined up to similarity for
$n=3$, and in particular when $r=3$.
\end{proof}

\begin{remark}
In the $r=2$ case, this observation was noted in \cite[Lemma
4.2]{GQMT01}.
\end{remark}

\begin{lem} \label{e_n.w(alpha)}
Let $\phi$ be an $r$-fold Pfister form, and $q$, $q'$ quadratic
forms over $k$. Then $\phi \otimes q \cong \phi \otimes q'$
implies $e_r(\phi)w(q)=e_r(\phi)w(q') \in H^*(k).$
\end{lem}

\begin{proof}
This is an extension of \cite{De62} or \cite{Mi70} where it is
shown for $r=0$. We need the fact from \cite{WS77} that says if
$\phi \otimes q \cong \phi \otimes q'$ then $q$ and $q'$ are
$\phi$-chain equivalent.

Two quadratic forms are \textit{simply $\phi$-equivalent} if they
can both be diagonalized in such a way that $q \cong \langle a_1,
\cdots, a_n \rangle$ and $q'' \cong \langle \lambda a_1, a_2,
\cdots , a_n \rangle$, where $\lambda$ is represented by $\phi$.
Then two forms are \textit{$\phi$-chain equivalent} if there is a
finite chain of simple $\phi$-equivalences from one to the other.

This immediately reduces the problem to showing that equality
holds at each stage of the chain equivalence. This is the same as
showing $e_r(\phi)w(\langle a \rangle)= e_r(\phi)w(\lambda\langle
a \rangle)$ for $\lambda$ represented by $\phi$. For such a
$\lambda$, we have $\phi \otimes \langle 1,-\lambda \rangle$ is
isotropic, and hence $e_r(\phi) \cdot (\lambda) =0 \in
H^{r+1}(k)$. Expanding $w(\lambda \langle a \rangle)=1+(\lambda) +
(a)$, the result clearly follows.
\end{proof}

So the following Lemma shows that the quadratic form $\phi \otimes
q$, where $\textup{det}(q)=1$, is in fact determined up to
isometry by the isomorphism class of a reduced Jordan algebra.
Since $q$ is odd dimensional, there is always a determinant 1
quadratic form similar to it. We will write $d(q)=\textup{det}(q)$
for the element of $k^*/(k^*)^2$ corresponding to $w_1(q)\in
H^1(k)$.

\begin{lem} \label{UptoIsometry}
Let $\phi$ be an $r$-Pfister form, $\lambda \in k^*$, and $q,q'$
quadratic forms. If $\phi \otimes q \cong \phi \otimes \lambda q'$
then $\phi \otimes d(q)q \cong \phi \otimes d(q')q'.$
\end{lem}

\begin{proof}
\begin{align*}
e_r(\phi) \cdot w_1(q) = e_r(\phi) \cdot w_1(\lambda q') \in
H^{r+1}(k) &\Leftrightarrow e_r(\phi) \cdot (\lambda d(q) d(q')) =
0
\\ &\Leftrightarrow \phi \otimes \langle \langle \lambda d(q)
d(q') \rangle \rangle \textup{ is hyperbolic} \\ & \Leftrightarrow
d(q)d(q')=\lambda \textup{ mod }D(\phi)^* \\ &\Rightarrow \phi
\otimes d(q)q \cong \phi \otimes d(q')q'.
\end{align*}
\end{proof}

\begin{proof}[Proof of Thm.\ \ref{e_r.w}]
First we will show that the invariants $e_r \otimes w_{2i}$ extend
to invariants on $k$-isomorphism classes of reduced Jordan
algebras.

Consider the reduced Jordan algebra $J=(M_n(C),ad_q \otimes -)_+$
with $n=2m+1$. Then we can assume $\textup{det}(q)=1$. Lemma
\ref{UptoSimilarity} together with Lemma \ref{UptoIsometry} show
that $\phi \otimes q$ is determined up to isometry by the
isomorphism class of $J$. Then by Lemma \ref{e_n.w(alpha)} we can
define $v_i(J)=e_r(\phi)w_{2i}(q)\in H^{r+2i}(k)$ on
$k$-isomorphism classes of reduced Jordan algebras, for $1\leq i
\leq m$. This clearly extends $e_r \otimes w_{2i}$.

Finally, by Lemma \ref{ReducingFields}, any odd degree Jordan
algebra becomes reduced after an odd degree extension. So by
\cite[Prop.\ 7.2]{Ga06} these invariants may be extended to
non-reduced Jordan algebras as well, and in other words, to all
$\mathbf{Aut}(J^r_n)$-torsors. By Cor.\ \ref{InvInject}, $v_i$ is
the unique invariant extending $e_r \otimes w_{2i}$.
\end{proof}

\begin{remark}
For $r=1$, there is an invariant closely related to $v_1$ defined
on conjugacy classes of algebras with unitary involution of odd
degree in \cite[p.438, 31.44]{KMRT}. They related it to the Rost
invariant.
\end{remark}

Now we can state and prove our main theorem.

\begin{thm} \label{AutInvs}
$\textup{Inv}_{k_0}(\mathbf{Aut}(J^r_n))$ is a free
$H^*(k_0)$-module with a basis consisting of the invariants $\{1,
v_0, v_1, v_2, \cdots , v_m \}$.
\end{thm}

\begin{proof}
For $r=0$ this is shown in \cite[ch.\ VI]{GMS03}, noting that in
this case $v_0=1$, causing a redundancy in the set of basis
elements. So take $r>0$. From Cor.\ \ref{InvInject} we know that
every $\mathbf{Aut}(J^r_n)$-invariant restricts to some $$1
\otimes a + e_r \otimes b \in \textup{Inv}(\mathsf{Pf}_r) \otimes
\textup{Inv}(\mathsf{Quad}_{n,1}),$$ for some uniquely defined
$a,b \in \textup{Inv}(\mathsf{Quad}_{n,1})$. We know from
\cite[ch.\ VI]{GMS03} that any $b \in
\textup{Inv}(\mathsf{Quad}_{n,1})$ is in the $H^*(k)$-span of the
even Stiefel-Whitney classes, so by Thm.\ \ref{e_r.w}, $e_r
\otimes b$ is the restriction of some
$\mathbf{Aut}(J^r_n)$-invariant in the $H^*(k)$-span of $\{v_0,
v_1, \cdots , v_m \}$. So all that remains to show is that if $1
\otimes a$ is the restriction of an
$\mathbf{Aut}(J^r_n)$-invariant, then $a$ is constant.

Let $a'$ be an $\mathbf{Aut}(J^r_n)$-invariant that restricts to
$1 \otimes a$ for some $\mathsf{Quad}_{n,1}$-invariant $a$. If we
let $C_s$ be the split composition algebra of dimension $2^r$,
then $$J=(M_n(C_s),ad_q \otimes -)_+$$ is isomorphic to the split
algebra $J^r_n$ (by Prop. \ref{ad_h-h} for $r\neq 3$ and by
\cite[Cor.\ 5.8.2]{SV00} for $r=3$). So for any such $J$, we must
have that $a'(J)=a(q)$ is a constant, independent of $q$. Since we
can take $q$ to be an arbitrary $n$-dimensional form of
determinant 1, this implies $a$ is constant. This completes the
proof.
\end{proof}

\begin{remark}
One may ask to what extent do these $v_i$ determine the
$\mathbf{Aut}(J^r_n)$-torsors? There are examples of non-isometric
quadratic forms with determinant 1 in each dimension $\geq 4$ that
have equal total Stiefel-Whitney classes \cite[Beispiel
3.4.1]{Sch67}. So one can use these examples to write down two
different reduced $\mathbf{Aut}(J^r_n)$-torsors for $n\geq 4$ odd,
whose invariants agree.

In the case of $n=3$, on the other hand, for $r=0$ or 2, the
torsors are determined completely by $v_0$ and $v_1$. This is
because they are determined by their quadratic trace form
\cite[\S5]{SV00}. But for $r=1$ and $r=3$ this is not the case,
because (for $n=3$) the trace form of any non-reduced algebra is
isometric to the trace form of some reduced algebra. Nevertheless,
for $r=1$ one may define a degree 2, mod 3 invariant, which
together with $v_0$ and $v_1$, classify
$\mathbf{Aut}(J^1_3)$-torsors \cite[\S19.B, \S30.C]{KMRT}. For
$r=3$, one may define a degree 3, mod 3 invariant, but it is an
open problem whether this invariant together with $v_0$ and $v_1$,
classify $\mathbf{Aut}(J^3_3)$-torsors \cite[9.4]{Se95}.
\end{remark}

\section{Essential dimension}

Given an algebraic group $G$ over $k_0$, and a $G$-torsor $E$ over
$k$, the essential dimension of $E$ is defined to be the minimum
transcendence degree over $k_0$ of all fields of definition of
$E$. The essential dimension of an algebraic group is defined to
be the supremum of the essential dimensions of all of its torsors
(\cite{RY00}, \cite{BF03}, \cite{CS06}). The essential dimension
of most simple algebraic groups is unknown. We will determine the
value of the essential dimension at $2$ for some groups $G$, which
we will denote by ed$(G;2)$ (see \cite{CS06} for a definition of
the essential dimension at a prime). In all cases that we
consider, ed$(G;2)$ is equal to the lower bound given in
\cite[Theorem 1]{CS06} (or \cite[Thm.\ 7.8]{RY00} for
characteristic 0).

\begin{prop} \label{essdim}
For $n \geq 3$ odd, we have \textup{ed}$(\mathbf{Aut}(J^r_n); 2)=
r+n-1$.
\end{prop}

\begin{proof}
By the surjectivity in Lemma \ref{Surject} we have that for any
$\mathbf{Aut}(J^r_n)$-torsor $J$ over $k$, there is an odd degree
extension $L/k$ such that $J_L$ is reduced. So by using
\cite[Lemma 1.11]{BF03} we have that ed$(J;2) \leq
\textup{ed}_L(J_L) \leq \textup{ed}(\mathsf{Pf_r}) +
\textup{ed}(\mathsf{Quad}_{n,1})=r+n-1$. This gives us the upper
bound ed$(\mathbf{Aut}(J^r_n); 2)\leq r+n-1$.

The lower bound follows from \cite[Theorem 1]{CS06}.
Alternatively, we could deduce the lower bound by using the
non-triviality of the degree $r+n-1$ cohomological invariant
$v_m$. This follows from a slight modification of \cite[Cor.\
3.6]{BF03}, that if there is a degree $d$ invariant mod $2$, then
the essential dimension at $2$ is at least $d$.
\end{proof}

Let us consider what Prop.\ \ref{essdim} says for different $r$
and $n\geq 3$ odd.

For $r=0$ we get the well known fact that ed$(\mathbf{SO}_n)=
\textup{ed}(\mathbf{SO}_n;2)=n-1$.

For $r=1$ we get ed$(\mathbb{Z}/2\mathbb{Z} \ltimes
\mathbf{PGL}_n)\geq \textup{ed}(\mathbb{Z}/2\mathbb{Z} \ltimes
\mathbf{PGL}_n;2)=n$. The exact value of
ed$(\mathbb{Z}/2\mathbb{Z} \ltimes \mathbf{PGL}_n)$ is unknown to
the author for any $n\geq 3$.

For $r=2$ we get ed$(\mathbf{PSp}_{2n})=
\textup{ed}(\mathbf{PSp}_{2n};2)=n+1$, since all
$\mathbf{PSp}_{2n}$-torsors are reduced. Previously, the best
known upper bound for ed$(\mathbf{PSp}_{2n})$ was $2n^2-3n-1$,
which holds for $n$ even as well \cite{Le04}.

For $r=3$ we get ed$(F_4)\geq \textup{ed}(F_4;2)=5$, which is the
best known lower bound for ed$(F_4)$. The best published upper
bound for the essential dimension is ed$(F_4)\leq 19$ in
\cite{Le04}. \cite{Ko00} claimed to show that ed$(F_4)=5$, but
there was a mistake in the proof.

\section{Acknowledgements}

I would like to thank Jean-Pierre Tignol for pointing out the
$r=2$ case of Lemma \ref{ReducedAlg}, and Skip Garibaldi for
suggesting the application to essential dimension, and catching
several mistakes in a previous version of this paper. I would also
like to thank my Ph.D.\ supervisor Burt Totaro, as well as Carl
McTague and Arthur Prendergast-Smith for many useful discussions.
I also thank the referee for their comments and suggestions.

DPMMS, University of Cambridge, Wilberforce Road, Cambridge,
CB\textup{3 0}WB, United Kingdom \\ e-mail\textup{:
\texttt{M.MacDonald@dpmms.cam.ac.uk}}

\end{document}